\documentclass[a4paper,twoside]{article}
\usepackage{amssymb,amsmath,amscd,amsthm}

\input{epsf}

\newcommand{\card}[1]{{\mid\! #1 \!\mid}}
\newcommand{\pro}[2]{\langle #1, #2 \rangle}

\def\lolra{\Longleftrightarrow}

\def\N{{\mathbb N}}
\def\Nspos{{\N_{>0}}}
\def\Nato{\Nspos}
\def\weights{\N_{>0}^{d+1}}
\def\P{{\mathbb P}}
\def\Q{{\mathbb Q}}
\def\R{{\mathbb R}}
\def\Z{{\mathbb Z}}

\def\Qred{{Q_{\rm red}}}
\def\Pred{{P_{\rm red}}}
\def\red{{\rm red}}
\def\lcm{{\rm lcm}}

\def\F{{\cal F}}

\def\V{{{\cal V}}}

\def\conv{{\rm conv}}

\def\pos{{\rm pos}}

\def\vol{{\rm Vol}}

\def\NR{{N_{\R}}}

\def\MR{{M_{\R}}}

\def\fan{\triangle}

\def\KX{{K_X}}

\newtheorem*{theoremCs*}{Theorem C'}
\newtheorem*{theoremBs*}{Theorem B'}
\newtheorem*{theoremAs*}{Theorem A'}
\newtheorem*{theoremC*}{Theorem C}
\newtheorem*{theoremB*}{Theorem B}
\newtheorem*{theoremA*}{Theorem A}
\newtheorem*{theorem*}{Theorem}
\newtheorem*{corollary*}{Corollary}
\newtheorem{theorem}{Theorem}[section]
\newtheorem{definition}[theorem]{Definition}

\newtheorem{example}[theorem]{Example}
\newtheorem{corollary}[theorem]{Corollary}
\newtheorem{proposition}[theorem]{Proposition}
\newtheorem{lemma}[theorem]{Lemma}
\newtheorem{conjecture}[theorem]{Conjecture}

\long\def\symbolfootnote[#1]#2{\begingroup%
\def\thefootnote{\fnsymbol{footnote}}\footnote[#1]{#2}\endgroup}

\title{Volume and lattice points of reflexive simplices}

\author{\vspace*{-1ex}\sc \normalsize{Benjamin Nill} \\[2ex]
\small  \em Research Group Lattice Polytopes, FU Berlin\\[-.5ex]
\small \em Arnimallee 3, 14195 Berlin, Germany  \\[-.5ex]
\small \em e-mail: nill@math.fu-berlin.de \\[-.5ex]
 }

\begin{document}

\date{}

\maketitle
\vspace*{-5ex}

\begin{abstract}
Using new number-theoretic bounds on the denominators of unit fractions summing up to one, 
we show that in any dimension $d \geq 4$ there is only one $d$-dimensional reflexive simplex having maximal volume. 
Moreover, only these reflexive simplices can admit an edge that has the maximal number of lattice points possible for an edge 
of a reflexive simplex. 
In general, these simplices are also expected to contain the largest number of lattice points even among all lattice polytopes 
with only one interior lattice point. 
Translated in algebro-geometric language, our main theorem yields 
a sharp upper bound on the anticanonical degree of $d$-dimensional $\Q$-factorial Gorenstein toric Fano varieties with Picard number one, 
e.g., of weighted projective spaces with Gorenstein singularities.
\end{abstract}

\section*{Introduction}
\label{intro}

To make this article as accessible as possible to readers from different backgrounds, we summarize our results 
from the points of view of convex geometry, algebraic geometry and number theory respectively.

\smallskip

\textbf{Convex and discrete geometry:} 
Given\symbolfootnote[0]{{\em Mathematics Subject Classification (2000)}: Primary 14M25; Secondary 11D75, 11H06, 52A43, 52B20, 52C07.} 
a $d$-dimensional lattice polytope $P$ with only one lattice point in its interior, 
Hensley \cite{Hen83} showed that the volume and the number of lattice points of $P$ is bounded above by 
a function depending only on $d$. In \cite{LZ91} and \cite{Pik01} asymptotically better 
upper bounds were obtained. However, in lower dimensions they are presumably still very far from being sharp, 
cf. \cite{Kas06} and \cite{Rez06}. Nevertheless, 
in any dimension $d \geq 3$ there exists an explicit candidate that is expected to yield 
the maximal values. This is a simplex, here called $S_{Q'_d}$, 
that was described already more than twenty years ago in \cite{ZPW82}: $S_{Q'_d}$ is (up to translation) the convex hull of the origin 
of the lattice and $y_0 e_0$, \ldots, $y_{d-2} e_{d-2}$, $2(y_{d-1}-1) e_{d-1}$, where $e_0, \ldots, e_{d-1}$ is a lattice basis and 
$(y_0, y_1, y_2, y_3\ldots) = (2,3,7,43, \ldots)$ the so-called Sylvester sequence. 

In \cite{Bat94}, Batyrev introduced reflexive polytopes, which are lattice polytopes $P$, having the origin of the lattice in its interior, 
whose dual polytopes $P^*$ are also lattice polytopes. 
From the convex-geometric point of view, reflexive polytopes are special lattice polytopes with only one interior lattice point, 
and, since there are already in dimension four nearly half a billion isomorphism classes \cite{KS00}, 
they can be seen as a good testing ground for more general conjectures. 

Recently, Haase and Melnikov noticed in \cite{HM04} that $S_{Q'_d}$ is actually a reflexive simplex. 
Now, in the present paper we prove that in dimension $d \geq 4$ the simplices $S_{Q'_d}$ are the only reflexive 
simplices having maximal volume (Thm. A). This yields also an upper bound on the number of lattice points of reflexive simplices (Cor. \ref{ptcoro}). 
Moreover, we show that there is a sharp bound on the number of lattice points of an edge of a reflexive simplex, which 
is only obtained by $S_{Q'_d}$ (Thm. B). 
This was first observed by Haase and Melnikov in \cite{HM04} for $d \leq 4$. Finally, we give a 
sharp upper bound on the product of the volumes of dual reflexive simplices (Thm. C). This can be seen as 
an isoperimetric inequality of Blaschke-Santal\'o type for reflexive simplices.

\smallskip

\textbf{Algebraic and toric geometry:} The Borisov-Alexeev conjecture states that $d$-dimensional Fano varieties with 
$\epsilon$-log terminal singularities form a bounded family, in particular their anticanonical degree is bounded. 
This conjecture is of importance in birational geometry, and was proven in the toric case in \cite{BB93}. However good upper bounds 
are known only in particular cases. For instance, it was conjectured that $(d+1)^d$ is the maximal degree of a 
nonsingular Fano variety with Picard number one, achieved only for $\P^d$. 
This was proven for $d \leq 4$ in \cite{Hwa03}, however it is wrong without the assumption of nonsingularity. For higher Picard numbers this bound 
does not hold anymore, see \cite{Bat82}. In \cite[Thm.9]{Deb03} Debarre showed as a good approximation 
$(-K_X)^d \leq d! d^{\rho d}$ for a $d$-dimensional nonsingular toric Fano variety $X$ with Picard number $\rho \geq 2$.

Now, in \cite{Pro05} an effective bound (namely $72$) was 
given for canonical Gorenstein Fano threefolds, proving the so-called Fano-Iskovskikh-conjecture. 
Proposing the conjectural value in higher dimensions, 
we provide here in Thm. A' a sharp upper bound (namely $2 (y_{d-1} - 1)^2$) on the anticanonical degree $(-K_X)^d$ of 
(necessarily canonical) Gorenstein toric Fano varieties $X$ of dimension $d \geq 4$, under the additional assumptions 
that $X$ is $\Q$-factorial with Picard number one. Moreover, we have an effective bound on the anticanonical degree $-K_X . C$ 
of torus-invariant integral curves $C$ on such $X$ (Thm. B'). In both theorems we can completely characterize the cases of equality. 

These results base on the correspondence of Gorenstein toric Fano varieties via toric geometry to reflexive polytopes. 
These objects have been thoroughly studied, since the natural duality of reflexive polytopes yields a construction of mirror-symmetric Calabi-Yau varieties 
as general anticanonical hypersurfaces of the associated Gorenstein toric Fano varieties, see \cite{Bat94}.

Stimulating further research we would like to pose here the following question: While 
Thm. A' determines in any dimension the maximal anticanonical degree of a weighted projective space with Gorenstein singularities, 
we would like to know what is the minimal 
anticanonical degree a $d$-dimensional weighted projective space with Gorenstein singularities can have? Could it be $(d+1)!$ ?

\smallskip

\textbf{Elementary number theory:} 
It is well-known that questions concerning the volume and lattice points of simplices 
are closely related to geometry of numbers and diophantine equations. Here, the crucial objects of investigation are 
finite families of positive natural numbers whose reciprocals sum up to one,  e.g., $(2,3,6)$, since $1/2 + 1/3 + 1/6 = 1$. 
We call these unit partitions. In the study of these objects the Sylvester sequence plays a crucial role, as 
exemplified by a famous result of Curtiss \cite{Cur22}. Now, the proof of above results is based on new inequalities on unit partitions 
(Thm. \ref{kprop}). They were achieved by a method of proof first used by Izboldin and Kurliandchik in \cite{IK95}.

\smallskip

\textbf{Organization of the paper:} 
This paper is organized in the following way:

In the first section we start with the basic definitions and the main convex-geometric results (Thms. A,B,C). 
In the second section we give their algebro-geometric formulation (Thms. A',B',C'). 
The third section contains a discussion of the correspondence between 
weighted projective spaces, weight systems and lattice simplices containing the origin in their interiors. 
The fourth section then relates reflexive simplices and unit partitions, i.e., unit fractions summing up to one. 
In the fifth section we present our main number-theoretic result (Thm. \ref{kprop}). 
Finally, the last section contains the proofs of Theorems A,B,C.

\smallskip

\textbf{Remark:} One should note that no knowledge of algebraic or toric geometry is necessary to understand the convex-geometric 
results and their proofs. Moreover, basic facts from toric geometry that are freely used in this article 
like the dictionary of fans/lattice polytopes and toric varieties can be found 
in \cite{Ful93}. All relevant properties of reflexive polytopes like their correspondence to 
Gorenstein toric Fano varieties are contained in \cite{Bat94} or \cite{Nil05}. 

\section{The main results, convex-geometric}

{\em Throughout the paper let $d \geq 2$.}

\smallskip

Let $N,M$ be dual lattices of rank $d$, with the pairing $\pro{\cdot}{\cdot} \,:\, M \times N \to \Z$. 
We define $\NR := N \otimes_\Z \R$ and $\MR := M \otimes_\Z \R$. 
A {\em lattice polytope} (in $\MR$) is a polytope 
whose vertices are lattice points (in $M$). Two lattice polytopes $P_1$ and $P_2$ are regarded as {\em isomorphic} or {\em unimodularly equivalent}, if 
there is a lattice isomorphism inducing an affine bijection between the affine hulls of $P_1$ and $P_2$ which maps 
every $i$-dimensional face of $P_1$ onto an $i$-dimensional face of $P_2$, for all $i$, $0 \leq i \leq \dim(P_1) = \dim(P_2)$. 
In this case we write $P_1 \cong P_2$. Moreover, if $P \subseteq \MR$ is a lattice polytope, then 
$\vol(P)$ denotes the (normalized) {\em volume} of $P$, i.e., $d!$ times the euclidean volume of $P$ with respect to a lattice basis of $M$.

\smallskip

In \cite{Bat94}, Batyrev introduced the following notion:

\begin{definition}{\rm 
A {\em reflexive polytope} $P \subseteq \MR$ is a $d$-dimensional lattice polytope containing the origin in its interior 
whose dual polytope $P^*$ is also a lattice polytope, where the dual polytope is defined as
$$P^* := \{x \in \NR \,:\, \pro{y}{x} \geq -1 \,\forall\, y \in P\}.$$
}
\end{definition}

Hence, duals of reflexive polytopes are again reflexive.

\smallskip

Now, we need an elementary number-theoretic notion in order to construct some special reflexive simplices:

\begin{definition}{\rm 
The well-known recursive sequence \cite[A000058]{Slo04} of pairwise coprime natural numbers 
$y_0 := 2$, $y_n := 1 + y_0 \cdots y_{n-1}$ ($n \geq 1$) is called {\em Sylvester sequence}. 
It satisfies $y_n = y_{n-1}^2 - y_{n-1} + 1$ and starts as 
$y_0 = 2$, $y_1 = 3$, $y_2 = 7$, $y_3 = 43$, $y_4 = 1807$. 
We also define $t_n := y_n - 1$ for $n \in \N_{\geq 0}$.
\label{sylvester}
}
\end{definition}

Using these numbers we define two special families:

\pagebreak

\begin{definition}{\rm 
\begin{itemize}
\item[]
\item The $d+1$-tuple of natural numbers 
$$Q_d := (\frac{t_d}{y_0}, \ldots, \frac{t_d}{y_{d-1}}, 1)$$ 
is called {\em Sylvester weight system} of length $d$.
\item The $d+1$-tuple of natural numbers 
$$Q'_d := (\frac{2 t_{d-1}}{y_0}, \ldots, \frac{2 t_{d-1}}{y_{d-2}}, 1, 1)$$ 
is called {\em enlarged Sylvester weight system} of length $d$.
\end{itemize}
\label{modisylv}
}
\end{definition}

We can define lattice simplices from these weight systems:

\begin{definition}{\rm 
Let $Q$ be a $d+1$ tuple $(q_0, \ldots, q_d)$ of positive natural numbers such that $q_d = 1$. Let $\card{Q} := q_0 + \cdots + q_d$. 
If $k_i := \card{Q}/q_i$ is a natural number for any $i = 0, \ldots, d-1$, then we define 
\[S_Q := \conv(k_0 e_0 - u, \ldots, k_{d-1} e_{d-1} - u, -u),\]
for $u := e_0 + \cdots + e_{d-1}$, where $e_0, \ldots, e_{d-1}$ is a $\Z$-basis of $M$.
\label{explicit}
}
\end{definition}

For such $Q$ (in particular for $Q_d$ and $Q'_d$) the associated simplices $S_Q$ are reflexive, see Prop. \ref{constr}. 
Now, here is our first main result:

\begin{theoremA*}
\begin{enumerate}
\item[]
\item $S_{(1,1,1)}$ is the unique two-dimensional reflexive simplex with the largest volume $9$, 
respectively the largest number of lattice points $10$. 
\item $S_{(3,1,1,1)}$ and $S_{Q'_3}$ are the only three-dimensional reflexive simplices with the largest 
volume $72$, respectively the largest number of lattice points $39$. 
\item Let $d \geq 4$. Then $S_{Q'_d}$ is the unique $d$-dimensional reflexive simplex with the largest volume, namely $2 t_{d-1}^2$. 
\end{enumerate}
\label{theAs}
\end{theoremA*}

This bound vastly improves on general upper bounds on the volume of lattice simplices containing only one lattice point 
in the interior, as given in \cite{LZ91} or \cite{Pik01}. 

\smallskip

From the following theorem in \cite{Bli14} we get also an upper bound on the number of lattice points:

\begin{lemma}[Blichfeldt]
Let $P \subseteq \MR$ be a $d$-dimensional lattice polytope. Then $\card{P \cap M} \leq d + \vol(P)$.
\end{lemma}

Now, using a lower bound on the number of lattice points in $S_{Q'_d}$ due to \cite{ZPW82} 
we see that there is still a gap to bridge:

\begin{corollary}
Let $J$ denote the maximal number of lattice points some $d$-dimensional reflexive simplex can have. Then 
we get for $d \geq 3$:
\[\frac{1}{3 (d-2)!} t_{d-1}^2 < J \leq d + 2 t_{d-1}^2 \in O(c^{2^{d+1}}),\]
where $c \approx 1.26408$ (see Lemma \ref{double}).
\label{ptcoro}
\end{corollary}

We can explicitly state a sharper conjecture, that was checked 
by the computer classification of Kreuzer and Skarke \cite{KS04,KS05} for four-dimensional reflexive polytopes:

\begin{conjecture}{\rm 
Let $d \geq 4$: The reflexive simplex $S_{Q'_d}$ has solely the largest number of lattice points among all $d$-dimensional lattice polytopes 
with only one interior lattice point.
}
\end{conjecture}

Our second main result is another indication in this direction. It has been observed by Haase and Melnikov in \cite{HM04} 
for $d \leq 4$ and initialized this research:

\begin{theoremB*}
The maximal number of lattice points on an edge of a reflexive simplex of dimension $d$ is 
$2 t_{d-1} + 1$, with equality attained only for $S_{Q'_d}$.
\label{theBs}
\end{theoremB*}

Here is the third main result (its first part can be seen as a variant of the Blaschke-Santal\'o inequality for 
centred compact convex bodies, cf. \cite[p.165]{Lut93}):

\begin{theoremC*}
Let $P$ be a $d$-dimensional reflexive simplex. Then
\[(d+1)^{d+1} \leq \vol(P) \, \vol(P^*) \leq t^2_d,\]
with equality of the first bound iff the sum of the vertices of $P$ equals $0$, and equality of the second bound 
iff $P \cong S_{Q_d}$ (in this case $P \cong P^*$). 

\smallskip

Furthermore, let the vertices of $P^*$ generate the lattice $N$. Then
\[\vol(P^*) \leq t_d,\]
with equality iff $P^* \cong S_{Q_d} \cong P$.
\label{theCs}
\end{theoremC*}

\smallskip

The proofs of Theorems A,B,C will be given in the last section.

\section{The main results, algebro-geometric}
\label{main}

In this section we reformulate Theorems A,B,C algebro-geometrically. 

\smallskip

For $Q \in \weights$ we denote by $\P(Q)$ the associated weighted projective space (see \cite[p.35]{Ful93})
with weights $Q$, e.g., $\P(1,\ldots, 1) = \P^d$. For any $d \geq 2$ the weighted projective spaces $\P(Q_d)$ and 
$\P(Q'_d)$, where $Q_d$ and $Q'_d$ were defined in Def. \ref{modisylv}, 
have Gorenstein singularities, i.e., the anticanonical divisor is an ample Cartier divisor. This follows for instance from 
Theorem \ref{reftheo} below.

\smallskip

Here is our first main algebro-geometric result:

\begin{theoremAs*}
Let $X$ be a $d$-dimensional $\Q$-factorial Gorenstein toric Fano variety with Picard number one.
\begin{enumerate}
\item If $d=2$, then
\[(-\KX)^2 \leq 9,\]
with equality iff $X \cong \P^2$. 
\item If $d=3$, then
\[(-\KX)^3 \leq 72,\]
with equality iff $X \cong \P(3,1,1,1)$ or $X \cong \P(Q'_3) = \P(6,4,1,1)$.
\item If $d \geq 4$, then
\[(-\KX)^d \leq 2 t_{d-1}^2,\]
with equality iff $X \cong \P(Q'_d)$.
\end{enumerate}
\label{theA}
\end{theoremAs*}

Since Gorenstein toric Fano varieties have canonical singularities, this result motivates the following more general conjecture:

\begin{conjecture}{\rm 
The results of Theorem A' hold for Gorenstein Fano varieties with canonical singularities.
}
\end{conjecture}

In the case of threefolds the bound in Theorem A' is the so-called Fano-Iskovskikh-conjecture. 
It has very recently been proven by Prokhorov \cite{Pro05}. 

\smallskip

Recall that the {\em class number} of a $d$-dimensional toric variety $X$, i.e., the rank of the group of Weil divisors modulo
linear equivalence, is equal to the number of one-dimensional cones of the associated fan minus $d$ (see \cite[p. 63]{Ful93}). $X$ is {\em $\Q$-factorial}, i.e., 
the Picard number equals the class number, if and only if all cones of the associated fan are simplices. Now, a theorem of Batyrev \cite[Thm.4.1.9]{Bat94} yields that  
$X$ is a $d$-dimensional Gorenstein toric Fano variety if and only if the the associated fan $\fan$ 
(in N) is spanned by the faces of a $d$-dimensional reflexive polytope. 
In this case, $X$ is $\Q$-factorial with Picard number one, or equivalently, $X$ has class number one, if and only if $\fan$ is spanned by the faces of a 
reflexive simplex. Then we denote its dual reflexive simplex by $P$. Here, the vertices of the simplex $P^*$ generate the lattice $N$ if 
and only if $X$ is isomorphic to a weighted projective space (e.g., see Theorem \ref{mainweight}). 

\begin{definition}{\rm 
If $X = \P(Q)$ (for $Q \in \weights$) is a weighted projective space with Gorenstein singularities, then we 
define $S_Q$ as the (up to unimodular equivalence unique) reflexive simplex $P \subseteq \MR$. 
\label{deffy}}
\end{definition}

This definition is consistent with Def. \ref{explicit} due to Thm. \ref{reftheo} and Prop. \ref{constr}. 

\smallskip

Now, Theorem A' follows directly from Theorem A due to the following formula (see \cite[p.111]{Ful93}):
\begin{equation}
(-\KX)^d = \vol(P).
\label{volf}
\end{equation}

\smallskip

The second algebro-geometric theorem is concerned with the maximal anticanonical degree of integral torus-invariant curves 
(in contrast, in Mori theory usually one is rather interested in their minimal degree, cf. \cite[Prop.2.1]{Lat96}) :

\begin{theoremBs*}
Let $X$ be a $d$-dimensional $\Q$-factorial Gorenstein toric Fano variety with Picard number one. 
Let $C$ be a torus-invariant integral curve on $X$. Then 
\[-\KX.C \leq 2 t_{d-1},\]
where equality implies $X \cong \P(Q'_d)$.
\label{theB}
\end{theoremBs*}

\smallskip

To see how to obtain this result from Theorem B, we associate to $X$, as before, a fan $\fan$ in $N$ and a reflexive simplex $P \subseteq \MR$. 
We remark that any torus-invariant integral curve $C$ is given by a {\em wall} $\rho$, 
i.e., a $(d-2)$-dimensional cone of $\fan$. Moreover, $\rho$ is obviously in correspondence with an edge $e$ of $P$. 
Hence, Theorem B' follows from the following observation (see \cite[Cor. 3.6]{Lat96}):
\begin{equation}
-\KX.C = \card{e \cap M} - 1.
\label{edgesf}
\end{equation}

\smallskip

For the last main algebro-geometric result we need the following definition:

\begin{definition}{\rm 
If $P \subseteq \MR$ is a reflexive polytope and $X$ is the Gorenstein toric Fano variety associated to the fan (in $N$) spanned by the faces of $P^*$, 
we define $X^*$ as the Gorenstein toric Fano variety that 
is associated to the fan (in $M$) spanned by the faces of $P$.
}
\end{definition}

\begin{theoremCs*}
Let $X$ be a $d$-dimensional $\Q$-factorial Gorenstein toric Fano variety with Picard number one. Then
\[(d+1)^{d+1} \leq (-K_X)^d \, (-K_{X^*})^d \leq t^2_d,\]
where equality of the lower bound is attained only if $X$ is a quotient of $\P^d$ by a finite abelian group, and 
equality of the upper bound is attained only for $X \cong \P(Q_d)$ (in this case $X \cong X^*$).

\smallskip
Furthermore, let $X$ be a weighted projective space with Gorenstein singularities. Then
\[(-K_{X^*})^d \leq t_d,\]
with equality iff $X \cong \P(Q_d) \cong X^*$.
\label{theC}
\end{theoremCs*}

One could conjecture that the first part of this result, as well as Theorem B', might also be true for arbitrary Gorenstein toric Fano varieties.

\smallskip

Theorem C' is a direct translation of Theorem C (use (\ref{volf}), and the 
basic results of the next section, notably Thm. \ref{maps} and equation (\ref{charrel})).

\section{Lattice simplices and weight systems}

The next two sections summarize and extend results of Batyrev in \cite[5.4,5.5]{Bat94} and Conrads in \cite{Con02}. 

\smallskip

Throughout let $P \subseteq \MR$ be a $d$-dimensional rational simplex containing the origin of the lattice 
in the interior. We define $\V(P)$ as the set of {\em vertices} of $P$ and 
$\F(P)$ as the set of {\em facets} of $P$. For $F \in \F(P)$ we let $\eta_F \in \NR$ denote 
the unique inner normal of $F$ with $\pro{\eta_F}{F} = -1$, thus, $\eta_F \in \V(P^*)$. 
In general, $\conv(V)$ denotes the convex hull of a set $V \subseteq \MR$. 

\begin{definition}{\rm 
\begin{itemize}
\item[]
\item Let $P = \conv(v_0, \ldots, v_d)$. 
We define $q_i := \card{\det(v_j \,:\, j=0, \ldots, d, \, j \not=i)}$ for $i = 0, \ldots, d$. 
The family $Q_P := (q_0, \ldots, q_d)$ is called the associated 
{\em weight system} of $P$ of {\em total weight} $\card{Q_P} := q_0 + \cdots + q_d$. 
If $P$ is a lattice polytope, we define the {\em factor} of $P$ as the index of the sublattice of $M$ generated by 
the vertices $v_0, \ldots, v_d$. 
\item We explicitly define a {\em weight system} $Q$ {\em of length} $d$ as a $d+1$-tuple of positive 
rational numbers; two weight systems are regarded as {\em isomorphic}, if they are equal up to permutation. 
The {\em total weight} $\card{Q}$ is the sum of the elements of $Q$. 

Let $Q$ consist of natural numbers. 
Then we define the {\em factor} of $Q$ as $\lambda_Q := \gcd(Q)$. 
In this case $Q$ is called {\em reduced}, if $\lambda_Q = 1$. 
We denote by $\Qred$ its {\em reduction} $Q/\lambda_Q$. $Q$ is called {\em normalized}, if 
$\gcd(Q\backslash\{i\}) = 1$ for $i = 0, \ldots, d$.
\end{itemize}
}
\end{definition}

If $P$ is a lattice simplex, we have by \cite[Lemma 2.4]{Con02}
\[\lambda_P = \lambda_{Q_P} \in \N_{> 0}.\]
Furthermore, in this case $(Q_P)_\red = (q'_0, \ldots, q'_n)$ is the unique reduced weight system satisfying 
\begin{equation}
\sum_{i=0}^d q'_i v_i = 0.
\label{charrel}
\end{equation}

There is the following theorem (see \cite[4.4-7]{Con02}, also \cite[Thm. 5.4.5]{Bat94}):

\begin{theorem}[Batyrev, Conrads]
For a reduced weight system $Q$ exists up to isomorphism a unique lattice simplex $P_Q \subseteq \MR$ with $Q_{P_Q} = Q$. 
The toric variety associated to the fan spanned by the faces of $P_Q$ is isomorphic to $\P(Q)$.

Furthermore, if $P$ is a lattice simplex with $(Q_P)_\red = Q$, then there exists (in a lattice basis) 
a matrix $H$ in Hermite normal form with determinant $\lambda_P$ such that $P \cong H P_Q$. 
In particular, any complete toric variety with class number one is the quotient of a weighted projective space 
by the action of a finite abelian group.
\label{maps}
\end{theorem}

Regarding normalized weight systems gives uniqueness (see \cite[3.7,3.8]{Con02}):

\begin{theorem}[Conrads]
There is a correspondence between isomorphism classes of
\begin{itemize}
\item lattice simplices $P$ that contain the origin and whose vertices are primitive lattice points that 
generate the lattice
\item normalized weight systems $Q$
\item weighted projective spaces
\end{itemize}
that is given by $P \mapsto Q_P$ and $Q \mapsto \P(Q)$.
\label{mainweight}
\end{theorem}

\smallskip

There is now an important invariant of a weight system, this generalizes \cite[Def. 5.4]{Con02} 
as can be seen from Prop. \ref{refprop}(2):

\begin{definition}{\rm 
Let $Q$ be a weight system of length $d$. Then we define
\[m_Q := \frac{\card{Q}^{d-1}}{q_0 \cdots q_d} \in \Q_{>0}.\]
}\label{mqdef}
\end{definition}

For the next result we need a straightforward lemma:

\begin{lemma}

\begin{enumerate}
\item[]
\end{enumerate}
$\begin{array}{ll}
&\det\begin{pmatrix}(n_1 - 1) & -1 & \cdots & & & \cdots & -1 \\
&\ddots&\ddots&\ddots&&&\\
-1 & \cdots & -1 & (n_i - 1) & -1 & \cdots & -1\\
&&&\ddots&\ddots&\ddots&\\
-1 & \cdots & & & \cdots & -1 & (n_d -1)
\end{pmatrix}\\
= & n_1 \cdots n_d  - \sum_{j=1}^d \prod_{i=1, i\not=j}^d n_i.
\end{array}$

\label{detformel*}
\end{lemma}

Dualizing changes the associated weight system only by a factor:

\begin{proposition}
\[Q_{P^*} = m_{Q_P} Q_P.\]
\label{formula}
\end{proposition}

\begin{proof}

Let $Q := Q_P$, $t := \card{Q}$, 
$\V(P) = \{v_0, \ldots, v_d\}$. For $i = 0, \ldots, d$ we denote by $F_i = \conv(v_j \,:\, j = 0, \ldots, d, \, j \not=i)$ the facet of $P$ not 
containing $v_i$. 
Fix a $\Z$-basis of $M$ and its dual basis of $N$. 
Now let $i \in \{0, \ldots, d\}$, and let $A_i$ be the matrix consisting of 
the coordinates of $v_j$ ($j=0, \ldots, d$, $j\not=i$) as rows, and let 
$B_i$  be the matrix consisting of 
the coordinates of $\eta_{F_j}$ ($j=0, \ldots, d$, $j\not=i$) as columns. 
Since $v_i = -\sum_{j\not=i} \frac{q_j}{q_i} v_j$, we get $\pro{\eta_{F_i}}{v_i} = 
\sum_{j\not=i} \frac{q_j}{q_i} = \frac{t-q_i}{q_i} = \frac{t}{q_i} - 1$. 

Without restriction we set $i=0$. Applying the previous lemma to $A_0 B_0$, we get 
$\det(A_0 B_0) = \frac{t}{q_1} \cdots \frac{t}{q_d} (1 - \sum_{j=1}^d \frac{q_j}{t}) =
\frac{t}{q_1} \cdots \frac{t}{q_d} (1 - \frac{t-q_0}{t})
= \frac{t^{d-1}}{q_1 \cdots q_d} q_0$. 
Therefore $\det(B_0) = \frac{\det(A_0 B_0)}{\det(A_0)} = \frac{t^{d-1}}{q_1 \cdots q_d} = m_Q q_0$. 

\end{proof}

\section{Reflexive simplices and unit partitions}

Let's consider the case of a reflexive simplex. Here, the following definitions turn out to be convenient 
(the notion of unit partitions is closely related to the much-studied subject of Egyptian fractions, where the denominators 
have to be pairwise different, cf. \cite[D11]{Guy81}):

\begin{definition}{\rm 
\begin{itemize}
\item[]
\item A weight system is called {\em reflexive}, if it is reduced and any weight is a divisor of the total weight. 
Especially it has to be normalized.
\item A family of positive natural numbers $(k_0, \ldots, k_d)$ is called a {\em unit partition of total weight} 
$t' := \lcm(k_0, \ldots, k_d)$, if $\sum_{i=0}^d 1/k_i = 1$.
\end{itemize}
}
\end{definition}

There is the following observation (essentially due to Batyrev in \cite[5.4]{Bat94}), that is straightforward to prove:

\begin{proposition}
In the notation of the previous definition we have that 
mapping $(q_0, \ldots, q_d)$ to $(\frac{t}{q_0}, \ldots, \frac{t}{q_d})$, respectively 
mapping $(k_0, \ldots, k_d)$ to $(\frac{t'}{k_0}, \ldots, \frac{t'}{k_d})$ yields a bijection between 
reflexive weight systems and unit partitions.
\end{proposition}

The notion of a reflexive weight system is motivated by the following result \cite[Prop. 5.1]{Con02} 
(partially \cite[Thm. 5.4.3]{Bat94}):

\begin{theorem}[Batyrev, Conrads]
Under the correspondence of Theorem \ref{mainweight} we get correspondences of isomorphism classes of
\begin{itemize}
\item reflexive simplices whose vertices generate the lattice
\item reflexive weight systems, respectively, unit partitions
\item weighted projective spaces with Gorenstein singularities
\end{itemize}
\label{reftheo}
\end{theorem}

Recall that, if $Q$ is a reflexive weight system, then the reflexive simplex $S_Q$, as defined in Def. \ref{deffy}, is isomorphic to $(P_Q)^*$. 
It is difficult to give in general an explicit description of $S_Q$, cf. \cite{Con02}. 
However, in the case in which the vertices of some facet of $P_Q$ form a lattice basis, we can conveniently describe 
$P_Q$ and $S_Q$, cf. Def. \ref{explicit} (proof left to the reader):

\begin{proposition}
Let $Q$ be a reflexive weight system with weight $q_d = 1$. Denote the corresponding unit partition 
by $(k_0, \ldots, k_d)$. Then we have:
\[P_Q \cong \conv(e_0, \ldots, e_{d-1}, -q_0 e_0 - \cdots - q_{d-1} e_{d-1}),\]
\[S_Q \cong (P_Q)^* \cong \conv(k_0 e_0 - u, \ldots, k_{d-1} e_{d-1} - u, -u),\]
for $u := e_0 + \cdots + e_{d-1}$, where $e_0, \ldots, e_{d-1}$ is a $\Z$-basis 
of $M$.
\label{constr}
\end{proposition}

There is the following generalization of \cite[5.3, 5.5]{Con02}:

\begin{proposition}
Let $P \subseteq \MR$ be a reflexive simplex with associated weight system $Q := Q_P = \lambda_P \Qred$. 
Then $\Qred$ is a reflexive weight system. Let $\Pred \subseteq \MR$ be the corresponding reflexive simplex 
and $(k_0, \ldots, k_d)$ the associated unit partition. 
We have the following results:

\begin{enumerate}
\item $(Q_{P^*})_\red = \Qred$, $\lambda_{P^*} = m_{Q_P} \lambda_P = \frac{m_\Qred}{\lambda_P}$.
\item \[m_\Qred = \frac{k_0 \cdots k_d}{\lcm(k_0, \ldots, k_d)^2} = \lambda_{(\Pred)^*} \in \Nato.\]
\item 
There exist injective lattice homomorphisms inducing bijections
\[\Pred \stackrel{\lambda_P}{\to} P \stackrel{\frac{m_\Qred}{\lambda_P}}{\to} (\Pred)^* \cong S_{Q_\red},\]
where the number on the arrow gives the integer determinant of the respective lattice homomorphism.
\item $\lambda_P \;\mid\; m_\Qred$. Furthermore
\[\lambda_P = 1 \;\lolra\; P \cong \Pred,\quad \lambda_P = m_\Qred \;\lolra\; P \cong (\Pred)^*.\]
\item \[\vol(P) = \card{Q_P} \leq \vol(S_{\Qred}) = \frac{k_0 \cdots k_d}{\lcm(k_0, \ldots, k_d)},\]
with equality iff $P \cong S_{\Qred}$.
\end{enumerate}
\label{refprop}
\end{proposition}

\begin{proof}

We can easily see that $\Qred$ is a reflexive weight system (or use \cite[Prop. 5.1]{Con02}). 
1. Follows from \ref{formula}, (\ref{charrel}) and Def. \ref{mqdef}.
2. Follows from the construction of $(k_0, \ldots, k_d)$, Def. \ref{mqdef} and 1. 
3. First apply \ref{maps} to $P$. Then apply \ref{maps} to the lattice simplex $P^*$, 
use 1. and dualize. 
4. from 3.  
5. The (normalized) volume of a $d$-dimensional simplex having the origin as a vertex is the 
absolute value of the determinant of the matrix formed by the vertices different from the origin. 
Hence $\vol(P) = \card{Q_P} = \lambda_P \card{\Qred} = \lambda_P \,\lcm(k_0, \ldots, k_d)$. 
Now use 2. and 3.

\end{proof}

The second point of the proposition was already proven by Batyrev in \cite[Cor. 5.5.4]{Bat94} when 
investigating fundamental groups.

\smallskip

Corresponding to the special reflexive weight systems defined in Def. \ref{modisylv} we have some important unit partitions:

\pagebreak

\begin{definition}{\em 
\begin{itemize}
\item[]
\item $(y_0, \ldots, y_{d-1}, t_d)$ is 
called {\em Sylvester partition} of length $d$. 
It is a unit partition of total weight $t_d = t_{d-1} y_{d-1}$ corresponding to the reflexive weight system $Q_d$. 
It corresponds to a {\em self-dual} reflexive simplex, since $m_{Q_d} = 1$.
\item $(y_0, \ldots, y_{d-2}, 2 t_{d-1}, 2 t_{d-1})$ is called {\em enlarged 
Sylvester partition} of length $d$. It corresponds to $Q'_d$, and we have $m_{Q'_d} = t_{d-1}$.
\end{itemize}
}
\end{definition}

As an illustration of the previous notions we classify all five isomorphismclasses of two-dimensional reflexive simplices:

\begin{example}{\rm 
Let $d=2$. We have three unit partitions:

\begin{enumerate}
\item $(3,3,3)$ corresponding to $Q := (1,1,1)$. Here, $m_Q = 3$. So we have $P_Q \cong \conv((1,0),(0,1),(-1,-1)$ 
(corresponding to $\P^2$) and $S_Q \cong (P_Q)^* \cong$\\$\conv((2,-1),(-1,2),(-1,-1))$ 
as the only reflexive simplices $P$ with\\$(Q_P)_\red = Q$ (due to \ref{refprop}(4)).
\item The Sylvester partition $(2,3,6)$ corresponding to $Q := Q_2 = (3,2,1)$. Here, $m_Q = 1$. So this yields 
the self-dual reflexive simplex $P_Q \cong S_Q \cong$\\$\conv((1,0),(0,1),(-3,-2))$ 
(corresponding to $\P(Q)$). 
\item The enlarged Sylvester partition $(2,4,4)$ corresponding to $Q := Q'_2 = (2,1,1)$. 
Here, $m_Q = 2$. So we have $P_Q \cong \conv((1,0),(0,1),(-2,-1)$ 
(corresponding to $\P(Q)$) and $S_Q \cong (P_Q)^* \cong \conv((1,-1),(-1,3),(-1,-1))$ 
as the only reflexive simplices $P$ with $(Q_P)_\red = Q$.
\end{enumerate}
\label{twoweights}
}
\end{example}

\section{The main results, number-theoretic}

In this section we are going to prove the following result on unit partitions:

\begin{theorem}
Let $(k_0, \ldots, k_d)$ be a unit partition.
\begin{enumerate}
\item \[\lcm(k_0, \ldots, k_d)^2 \leq k_0 \cdots k_d \leq t_d^2,\]
with equality in the second case only for the Sylvester partition.

\item \[(d+1)^{d+1} \leq k_0 \cdots k_d,\]
with equality iff $(k_0, \ldots, k_d) = (d+1, \ldots, d+1)$.

\item Let $d \geq 3$ and $k_0 \leq \ldots \leq k_d$. Then
\[\frac{k_0 \cdots k_d}{\lcm(k_0, \ldots, k_d)} \leq k_0 \cdots k_{d-1} \leq 2 t_{d-1}^2,\] 
with equality in the second case iff $(k_0, \ldots, k_d)$ is the enlarged Sylvester partition or $(2,6,6,6)$. 
\end{enumerate}
\label{kprop}
\end{theorem}

First let's recall the most important classical result about unit partitions due to Curtiss \cite[Thm. I]{Cur22}:

\begin{theorem}[Curtiss]
Let $a_1, \ldots, a_m$ be positive integers such that $s := \sum_{i=1}^m \frac{1}{a_i}< 1$. Then
\[s \leq \sum_{i=0}^{m-1} \frac{1}{y_i} = 1 - \frac{1}{t_m},\]
with equality iff $\{a_1, \ldots, a_m\} = \{y_0, \ldots, y_{m-1}\}$.
\label{erdoes}
\end{theorem}

As a corollary we get that $\max(k_0, \ldots, k_d) \leq t_d$ for any unit partition $(k_0, \ldots, k_d)$, and 
equality holds only for the Sylvester partition. 

While the proof of Curtiss was quite complicated, the following nice result \cite[Lemma 1]{IK95} has a 
much simpler proof, and it easily yields again the previous theorem. 
Here we have included some statements that are implicit in their proof. 

\begin{lemma}[Izhboldin, Kurliandchik]
Let $x_1, \ldots, x_n$ be real numbers satisfying $x_1 \geq x_2 \geq \cdots \geq x_n \geq 0$, 
$x_1 + \cdots + x_n = 1$ and $x_1 \cdots x_k \leq x_{k+1} + \cdots + x_n$ for $k = 1, \ldots, n-1$. Then
\[x_n \geq \frac{1}{t_{n-1}},\; \quad x_1 \cdots x_n \geq \frac{1}{t_{n-1}^2},\]
where equality in the first case holds iff equality in the second case holds 
iff $x_i = \frac{1}{y_{i-1}}$ for $i = 1, \ldots, n-1$.
\label{russen}
\end{lemma}

Any unit partition fulfills (after reordering) the assumptions of the lemma: To see this let $(k_0, \ldots, k_d)$ be a unit partition with 
$k_0 \leq \ldots \leq k_d$. We define $x_i := 1/k_i$ for $i = 0, \ldots, d.$ 
Then for $i = 0, \ldots, d-1$ we get $0 < x_{i+1} + \cdots + x_d = 1 - x_0 - \cdots - x_i = 
\frac{k_0 \cdots k_i - \sum_{j=0}^i \prod_{l=0, l \not=j}^i k_l}{k_0 \cdots k_i} \geq \frac{1}{k_0 \cdots k_i} 
= x_0 \cdots x_i$.

\smallskip

Hence, the previous lemma proves the upper bound (and its equality case) in statement 1. of Theorem \ref{kprop}. 
The middle inequality in 1. follows from Prop. \ref{refprop}(2). Statement 2. of Theorem \ref{kprop} follows immediately from applying 
the inequality of arithmetic and geometric means to $1/k_0, \ldots, 1/k_d$. 

\smallskip

For the proof of the last point in Theorem \ref{kprop} we need a preliminary result:

\begin{lemma}
Let $n \geq 4$, and $1 \leq r \leq n-1$. Then
\[(r+1)^r \, t_{n-r-1}^{r+1} \leq 2 t_{n-2}^2,\]
with equality iff $r=1$ or $(n,r)=(4,2)$.
\label{cruc}
\end{lemma}

\begin{proof}

Proof by induction on $n$. By explicitly checking $n=4,5$, we can assume $n \geq 6$. 

For $r=1$ the statement is trivial, so let $r \geq 2$. 
By induction hypothesis for $(n-1,r-1)$ we have $r^{r-1} \, t_{n-r-1}^r \leq 2 t_{n-3}^2$, this yields 
$(r+1)^r \, t_{n-r-1}^{r+1} \leq 2 t_{n-3}^2 \, t_{n-r-1} (\frac{r+1}{r})^r r$. Since 
$(\frac{r+1}{r})^r < e$, it is enough to show $t_{n-3}^2 \, t_{n-r-1} \, e \, r \leq t_{n-2}^2$, 
or equivalently, $t_{n-r-1} \, e \, r \leq y_{n-3}^2$.

For $n \geq 6$ it is easy to see that $e (n-1) \leq y_{n-3}$ (e.g., by Lemma \ref{double} below). 
Hence, $t_{n-r-1} e r \leq t_{n-3} e (n-1) \leq t_{n-3} y_{n-3} < y_{n-3}^2$.

\end{proof}

The following lemma gives the 
doubly exponential behavior of the Sylvester sequence (e.g., \cite[(4.17)]{GKP89}):

\begin{lemma}
There is a constant $c \approx 1.2640847353\cdots$ (called {\em Vardi} constant, see 
\cite[A076393]{Slo04}) such that for any $n \in \N$
\[y_n = \left\lfloor c^{2^{n+1}} + \frac{1}{2} \right\rfloor.\]
\label{double}
\end{lemma}

Now using Lemma \ref{cruc} and the ideas of the proof of Lemma \ref{russen} we can show the following result, which 
immediately implies statement 3. in Theorem \ref{kprop}:

\begin{lemma}
Let $n \geq 4$, and $x := (x_1, \ldots, x_n)$ as in Lemma \ref{russen}. Then
\[x_1 \cdots x_{n-1} \geq \frac{1}{2 t_{n-2}^2},\]
with equality iff {\em either} $n=4$ and $x_1 = \frac{1}{2}$, $x_2 = x_3 = x_4 = \frac{1}{6}$, {\em or} $(1/x_1, \ldots, 1/x_n)$ equals 
the enlarged Sylvester partition of length $n-1$.
\label{meins}
\end{lemma}

\begin{proof} 

The proof splits into three parts.\\

\textbf{Part I:} 

Let $A$ denote the set of $n$-tupels $x$ satisfying the conditions of Lemma \ref{russen}. 
It is easy to see that we have for $x \in A$ necessarily $1 > x_1$ and $x_n > 0$. 

Since $A$ is compact, there exists some $x \in A$ with $x_1 \cdots x_{n-1}$ minimal. Because of 
$(\frac{1}{y_0}, \ldots, \frac{1}{y_{n-3}}, \frac{1}{2 t_{n-2}}, \frac{1}{2 t_{n-2}}) \in A$ we have 
$x_1 \cdots x_{n-1} \leq \frac{1}{2 t_{n-2}^2}$.

\smallskip

{\em We are going to show that $x_{n-1} = x_n$.}

\smallskip

So assume $x_{n-1} > x_n$.

\smallskip

{\em Claim:} In this case we have
\begin{eqnarray*}
x_1 > x_2 > \cdots > x_{n-1} > x_n,\\
x_1 \cdots x_k = x_{k+1} + \cdots + x_n \text{ for } k = 1, \ldots, n-2.
\end{eqnarray*}
{\em Proof of claim:}

By convention we set $x_0 := 1$ and $x_{n+1} := 0$. We proceed by induction on $l = 1, \ldots, n-1$, and 
assume by induction hypothesis that $x_1 > x_2 > \cdots > x_l \geq x_{l+1}$ and 
$x_1 \cdots x_k = x_{k+1} + \cdots + x_n$ for $k = 1, \ldots, l-2$. 

We distinguish three cases:

\begin{enumerate}
\item $x_l > x_{l+1}$ and $x_1 \cdots x_{l-1} = x_l + \cdots + x_n$. 

In this case we can proceed.

\item $x_l > x_{l+1}$ and $x_1 \cdots x_{l-1} < x_l + \cdots + x_n$.

This implies $l \geq 2$. Then we can find some $\delta > 0$ s.t. 
$x' \in A$ with $x'_{l-1} := x_{l-1} + \delta$, $x'_l := x_l - \delta$ and 
$x'_j := x_j$ for $j \in \{1, \ldots, n\}\backslash\{l-1,l\}$. Hence, 
$x'_1 \cdots x'_{n-1} = \frac{x_1 \cdots x_{n-1}}{x_{l-1} x_l} (x_{l-1} x_l + \delta (x_l - x_{l-1}) - \delta^2) 
< x_1 \cdots x_{n-1}$, a contradiction.

\item $x_l = x_{l+1} = \cdots = x_i > x_{i+1}$ for $l+1 \leq i \leq n-1$. 

This implies $l \leq n-2$. Again we find some $\delta > 0$ s.t. 
$x' \in A$ with $x'_l := x_l + \delta$, $x'_i := x_i - \delta$ and 
$x'_j := x_j$ for $j \in \{1, \ldots, n\}\backslash\{l,i\}$. 

This can be done, since otherwise there has to exist 
$l \leq j < i$ such that $x_1 \cdots x_j = x_{j+1} + x_{j+2} + \cdots + x_n$. Since $x_j = x_{j+1}$ we have 
$0 = (1 - x_1 \cdots x_{j-1}) x_j + x_{j+2} + \cdots + x_n$, a contradiction.

Since again $x'_1 \cdots x'_{n-1} < x_1 \cdots x_{n-1}$, we get a contradiction.

\end{enumerate}

{\em End of proof of claim.}
\medskip

So we have $x_1 = x_2 + \cdots + x_n = 1 - x_1$, hence $x_1 = \frac{1}{2} = \frac{1}{y_0}$. In the same way we have for $k = 2, \ldots, n-2$ 
that $x_1 \cdots x_k = x_{k+1} + \cdots + x_n = 1 - x_1 - \cdots - x_k$, hence 
we get by induction $\frac{1}{t_{k-1}} x_k = 1 - \frac{1}{y_0} - \cdots - \frac{1}{y_{k-2}} - x_k = \frac{1}{t_{k-1}} - x_k$, 
thus $x_k =  1 - t_{k-1} x_k$, this implies $x_k = \frac{1}{1 + t_{k-1}} = \frac{1}{y_{k-1}}$. So this yields
$x_1 = \frac{1}{y_0},\; \ldots,\; x_{n-2} = \frac{1}{y_{n-3}}$. 

Furthermore, $x_{n-1} + x_n = 1 - x_1 - \cdots - x_{n-2} = \frac{1}{t_{n-2}}$. Since $x_{n-1} > x_n$, we get 
$x_{n-1} > \frac{1}{2 t_{n-2}}$. Therefore we have proven 
$x_1 \cdots x_{n-1} > \frac{1}{2 t_{n-2}^2}$, a contradiction.

\medskip

Hence, $x_{n-1} = x_n$.\\

\textbf{Part II:}

Let $A'$ denote the set of $(n-1)$-tupels $w \in \R^{n-1}$ satisfying the following conditions: 
$w_1 \geq w_2 \geq \cdots \geq w_{n-2} \geq \frac{w_{n-1}}{2} \geq 0$, 
$w_1 + \cdots + w_{n-1} = 1$ and $w_1 \cdots w_k \leq w_{k+1} + \cdots + w_{n-1}$ for $k = 1, \ldots, n-2$.

Let $w \in A'$ be fixed with $w_1 \cdots w_{n-1}$ minimal. 

Since $(\frac{1}{y_0}, \ldots, \frac{1}{y_{n-3}}, \frac{1}{t_{n-2}}) \in A'$ due to $n \geq 4$, 
we have $w_1 \cdots w_{n-1} \leq \frac{1}{t_{n-2}^2}$.

\smallskip

{\em We are going to show $w = (1/2,1/6,1/3)$ or $w = (1/y_0, \ldots, 1/y_{n-3}, 1/t_{n-2})$; 
in particular, $w_1 \cdots w_{n-1} = \frac{1}{t_{n-2}^2}$.}

\smallskip

Let $z := w_s = \cdots = w_{n-2} = \frac{w_{n-1}}{2}$ for $1 \leq s \leq n-1$ minimal. 

We define $r := n-s$. There are three cases to consider:

\begin{enumerate}
\item $s=1$, i.e., $r=n-1$.

Then $1 = w_1 + \cdots + w_{n-2} + w_{n-1} = (n-2) z + 2 z = n z$, so $z = 1/n$. This implies 
$w_1 \cdots w_{n-1} = \frac{2}{n^{n-1}}$. However $\frac{1}{t_{n-2}^2} < \frac{2}{n^{n-1}}$ for $n \geq 4$ 
by \ref{cruc}, a contradiction.

\item $s=2$, i.e., $r=n-2$.

Then $1 = w_1 + w_2 + \cdots + w_{n-2} + w_{n-1} = w_1 + (n-3) z + 2z = w_1 + (n-1) z$, hence 
$w_1 = 1 - (n-1) z$. Since $w_1 > z$, we get $z < \frac{1}{n}$. On the other hand 
$w_1 \leq w_2 + \cdots + w_{n-1} = (n-1)z$, hence $z \geq \frac{1}{2(n-1)}$. 

We have $w_1 \cdots w_{n-1} = (1-(n-1) z) 2 z^{n-2}$. This function attains its minimum on the interval 
$[\frac{1}{2(n-1)}, \frac{1}{n}[$ only for $z = \frac{1}{2(n-1)}$. (The elementary analytical proof of this statement is left to the reader.)

Therefore 
$\frac{1}{t_{n-2}^2} \geq w_1 \cdots w_{n-1} \geq \frac{1}{{(2(n-1))}^{n-2}}$. However, from Lemma \ref{cruc} 
we get $\frac{1}{{(2(n-1))}^{n-2}} \geq \frac{1}{t_{n-2}^2}$, with equality only for $n=4$ and $r=2$. Hence this yields 
$z = \frac{1}{2(n-1)}$ and  $w = (\frac{1}{2},\frac{1}{6},\frac{1}{3})$.

\pagebreak
\item $s \geq 3$, i.e., $r \leq n-3$.

Now a similar reasoning as in the proof of the claim in Part I yields $$w_1 = \frac{1}{y_0},\; \ldots,\; w_{s-2} = \frac{1}{y_{s-3}}.$$ 

Then $1 = w_1 + \cdots + w_{s-2} + w_{s-1} + w_s + \cdots + w_{n-1} = 
1 - \frac{1}{t_{s-2}} + w_{s-1} + (n-s+1) z$, hence $w_{s-1} = \frac{1}{t_{s-2}} - (r+1) z$. 

Since $w_{s-1} > z$, we get $z < \frac{1}{(r+2) t_{s-2}}$.

Since $w_1 \cdots w_{s-2} w_{s-1} \leq w_s + \cdots + w_{n-1}$, we get 
$\frac{1}{t_{s-2}} (\frac{1}{t_{s-2}} - (r+1) z) \leq (r+1) z$, hence 
$z \geq \frac{1}{(r+1) t_{s-2}^2 (1 + \frac{1}{t_{s-2}})} = \frac{1}{(r+1) t_{s-1}}$.

Now we define $f(z) := \frac{1}{t_{s-2}} (\frac{1}{t_{s-2}} - (r+1) z) 2 z^r = w_1 \cdots w_{n-1}$. 

Since the function $f(z)$ is for $z > 0$ 
strictly monotone increasing up 
to some value and then strictly monotone decreasing, we see that\\
$\frac{1}{t_{n-2}^2} \geq w_1 \cdots w_{n-1} 
\geq \min(f(\frac{1}{(r+1) t_{s-1}}),f(\frac{1}{(r+2) t_{s-2}})) = $\\$\min(\frac{2}{(r+1)^r t_{s-1}^{r+1}}, 
\frac{2}{(r+2)^{r+1} t_{s-2}^{r+2}})$.

There are two cases:

\begin{enumerate}
\item $s = n-1$, i.e., $r=1$.

Here, 
$w_1 \cdots w_{n-1} \geq \min(\frac{1}{t_{n-2}^2}, \frac{2}{9 t_{n-3}^3})$. 
Due to $n \geq 4$ and $t_{n-2} = t_{n-3} (t_{n-3}+1)$ we have $\frac{1}{t_{n-2}^2} \leq \frac{2}{9 t_{n-3}^3}$, 
thus $w_1 \cdots w_{n-1} = \frac{1}{t_{n-2}^2}$, and 
$z = \frac{1}{2 t_{n-2}}$, since $z < \frac{1}{(r+2)t_{s-2}}$. Hence 
$w = (\frac{1}{y_0}, \ldots, \frac{1}{y_{n-3}}, \frac{1}{t_{n-2}})$.

\item $s \leq n-2$, i.e., $r \geq 2$.

Here, $\frac{2}{(r+1)^r t_{s-1}^{r+1}} \leq \frac{2}{(r+2)^{r+1} t_{s-2}^{r+2}}$ 
if and only if $\frac{y_{s-2}^{r+1}}{t_{s-2}} \geq \frac{(r+2)^{r+1}}{(r+1)^r}$. 
This is true for $r=2$. For $r \geq 3$ we have $(r+2)(\frac{r+2}{r+1})^r < (r+2) e < 3^r \leq 
y_{s-2}^r < \frac{y_{s-2}^{r+1}}{t_{s-2}}$. 
Hence $\frac{1}{t_{n-2}^2} \geq w_1 \cdots w_{n-1} \geq \frac{2}{(r+1)^r t_{s-1}^{r+1}}$, 
a contradiction to \ref{cruc}.

\end{enumerate}

\end{enumerate}

\textbf{Part III:}

{\em We are going to finish the proof.} 

\smallskip

Let $n \geq 4$ and $x = (x_1, \ldots, x_n)$ satisfy the conditions in Lemma \ref{russen}, where $x_1 \cdots x_{n-1}$ has minimal value. 
By Part I we have $x_{n-1} = x_n$. Hence, $w' := (x_1, \ldots, x_{n-2}, 2 x_{n-1})$ is an element in $A'$, so by Part II we get 
\[x_1 \cdots x_{n-1} = \frac{w'_1 \cdots w'_{n-1}}{2} \geq \frac{1}{2 t_{n-2}^2},\]
where equality yields $(x_1, \ldots, x_n) = 
(1/y_0, \ldots, 1/y_{n-3}, 1/(2 t_{n-2}), 1/(2 t_{n-2}))$ or 
$n=4$ and $(x_1,x_2,x_3,x_4) = (1/2,1/6,1/6,1/6)$.

\end{proof}

\section{Proofs of Theorems A,B,C}

Combining the results of the previous two sections we get:

\begin{corollary}
\begin{enumerate}
\item[]
\item If $Q$ is a reflexive weight system of length $d$, then
\[\card{Q} \leq t_d, \text{ with equality iff $Q$ is isomorphic to $Q_d$}.\]

\item If $P$ is a reflexive simplex, then
\[(d+1)^{d+1} \leq \card{Q_P} \card{Q_{P^*}} = \vol(P) \vol(P^*) \leq t^2_d,\]
with equality of the first bound iff the sum of the vertices of $P$ equals $0$, 
and equality of the second bound iff $P \cong S_{Q_d} (\cong S^*_{Q_d})$.

\item If $P$ is a reflexive simplex for $d\geq 3$, then
\[\card{Q_P} = \vol(P) \leq 2 t_{d-1}^2,\]
with equality iff $P \cong S_{\Qred}$ for $\Qred$ isomorphic to $Q'_d$ or $(3,1,1,1)$.
\end{enumerate}
\label{qbound}
\end{corollary}

\begin{proof}

1. Follows from Thm. \ref{kprop}(1), since $\card{Q} = \lcm(k_0, \ldots, k_d)$ for 
the corresponding unit partition. 2. We observe that for $Q_P = \lambda_P \Qred$ 
Prop. \ref{refprop}(1,2) yields $Q_{P^*} = \frac{m_\Qred}{\lambda_P} \Qred$ and $\card{Q_P} \card{Q_{P^*}} = k_0 \cdots k_d$, where $(k_0, \ldots, k_d)$ 
is the unit partition corresponding to $\Qred$. Now the statements follow from Prop. \ref{refprop}(5), Thm. \ref{kprop}(1,2), equation (\ref{charrel}), and the fact that 
$P_{Q_d} \cong S_{Q_d}$ is the only reflexive simplex $P$ with $(Q_P)_\red = Q_d$, since $m_{Q_d} = 1$. 
3. Follows from Prop. \ref{refprop}(5) and Thm. \ref{kprop}(3). 

\end{proof}

Eventually, we can finish the proofs of the three main theorems:

\begin{proof}[Proof of Theorem A]

The case $d=2$ was considered in Example \ref{twoweights}. For $d \geq 3$ use Corollary \ref{qbound}(3). 
Furthermore, note that for $d=3$ we have the well-known formula $ \card{P \cap M} = \vol(P)/2 + 3$, that can be 
derived for instance from \cite{Hib92}.

\end{proof}

In the proof of Theorem B we use the following observation (see \cite[2.2]{Fuj03}):

\begin{proposition}
Let $Q$ be a normalized weight system with associated simplex 
$P_Q = \conv(v_0, \ldots, v_d)$. 
Let $i,j \in \{0, \ldots, d\}$, $i \not= j$. 
Then we have for 
the torus-invariant integral curve $C$ on $\P(Q)$ associated to the wall $\pos(v_k \,:\, i \not= k \not=j)$:
\[(-K_{\P(Q)}).C = \frac{\card{Q}}{\lcm(q_i,q_j)}.\]

In particular, when $Q$ is a reflexive weight system with corresponding unit partition $(k_0, \ldots, k_d)$, 
we get:
\[(-K_{\P(Q)}).C = \gcd(k_i,k_j).\]
\label{fujino}
\end{proposition}

\begin{proof}[Proof of Theorem B]

The fact that $S_{Q'_d}$ satisfies the bound is trivial from the definition, see Prop. \ref{constr}. 

Let $P = \conv(v_0, \ldots, v_d) 
\subseteq \MR$ be a reflexive simplex with associated weight system $Q_P$, and let $Q := (Q_P)_{{\rm red}} 
= (q_0, \ldots, q_d)$ correspond to the unit partition $(k_0, \ldots, k_d)$. Let $v_d,v_{d-1}$ be the vertices of an edge $e$ of $P$ 
that contains the maximal number $\card{e \cap M}$. Thus, $c := \card{e \cap M} - 1 = \gcd(k_d,k_{d-1})$ by the formula (\ref{edgesf}) 
and Prop. \ref{fujino}. We may assume $k_d \geq k_{d-1}$ and $c \geq 2 t_{d-1}$.

Now Theorem \ref{erdoes} implies 
$1 - 1/k_{d-1} - 1/k_d = \sum_{i=0}^{d-2} 1/k_i \leq 1 - 1/t_{d-1}$, so 
$2/k_{d-1} \geq 1/k_{d-1} + 1/k_d \geq 1/t_{d-1}$. This implies $2 t_{d-1} \leq c \leq k_{d-1} \leq 2 t_{d-1}$, so we have $c = 2 t_{d-1} =  k_{d-1} = k_d$. 
In particular, we get $\sum_{i=0}^{d-2} 1/k_i = 1 - 1/t_{d-1}$, thus, Theorem \ref{erdoes} 
yields $\{k_0, \ldots, k_{d-2}\} = \{y_0, \ldots, y_{d-2}\}$. Hence $Q = Q'_d$.

\smallskip

By \ref{constr} we can choose 
$\V(\Pred) = \{e_0, \ldots, e_{d-1}, e_d := -q_0 e_0 - \cdots -q_{d-1} e_{d-1}\}$ for some 
$\Z$-basis $e_0, \ldots, e_{d-1}$ of $M$. By Theorem \ref{maps} there is 
(up to unimodular equivalence) 
a $d \times d$-matrix $H = \{h_{i,j}\}$ in Hermite normal form of determinant $\lambda_P$ 
such that $H(e_0) = v_0,  \ldots, H(e_{d-1}) = v_{d-1}, 
H(e_d) = v_d$. 
Recall that a quadratic matrix $H$ is in Hermite normal form, if it is a lower triangular matrix with natural numbers as 
coefficients such that $h_{i,j} < h_{j,j}$ for $i > j$.

Since $v_{d-1}$ is primitive, we get $h_{d,d} = 1$, so $v_{d-1} = e_{d-1}$. 
Because of 
$$\frac{1}{2 t_{d-1}} (H(e_d) - e_{d-1}) \in M \text{ and } \frac{q_i}{2 t_{d-1}} = 
\frac{1}{y_i} \text{ for }i = 0, \ldots, d-2$$ 
we get the following $d-1$ equations:
\begin{eqnarray*}
\frac{h_{1,1}}{y_0} \in \N\\
y_0 (\frac{h_{2,1}}{y_0} + \frac{h_{2,2}}{y_1}) \in \N\\
\ldots\\
y_0 \cdots y_{d-3} (\frac{h_{d-1,1}}{y_0} + \cdots + \frac{h_{d-1,d-1}}{y_{d-2}}) \in \N.
\end{eqnarray*}
Using the fact that $\gcd(y_i,y_j) = 1$ for $i \not= j$, we deduce by induction that 
$y_{i-1}$ divides $h_{i,i}$ for $i = 1, \ldots, d-1$. Hence, we have
\[m_Q = t_{d-1} = y_0 \cdots y_{d-2} \leq h_{1,1} \cdots h_{d-1,d-1} h_{d,d} = \det H = \lambda_P \leq m_Q.\] 
Now Prop. \ref{refprop}(4) implies $P \cong (\Pred)^* \cong S_{Q'_d}$.

\end{proof}

\begin{proof}[Proof of Theorem C]

The first statement follows from Corollary \ref{qbound}(2), and the second statement from Corollary \ref{qbound}(1) and $\vol(P) = \card{Q_P}$, since 
in this case $Q_P$ is already reduced, thus a reflexive weight system. Also note $P_{Q_d} \cong S_{Q_d}$.

\end{proof}

\subsection*{Acknowledgements}
The author would like to thank his thesis advisor Victor Batyrev for giving reference to \cite{Pro05}, 
C. Haase for his interest, D. Ploog for discussion, and A. Hermann for helpful conversation, as well as the anonymous referees and G. Ziegler for 
suggestions on the organization of the paper.

Most of this work is part of the author's thesis.

\end{document}